\documentclass[12pt]{article}
\def\today{September 15, 2000}
\usepackage{amssymb}
\usepackage{amsthm}
\newtheorem{theorem}{Theorem}[section]

\newtheorem{pr}[theorem]{Proposition}
\newtheorem{prob}[theorem]{Problem}

\theoremstyle{definition}
\newtheorem{de}[theorem]{Definition}

\theoremstyle{remark}

\newcommand{\Section}[1]{\section{#1}\setcounter{equation}{0}}

\newcommand{\F}{\mathbb{F}}
\newcommand{\C}{\mathbb{C}}

\newcommand{\tr}{\mathrm{tr}}

\newcommand{\mat}{M\! at_{n\times n}}

\newcommand{\cZ}{\mathcal{Z}}
\newcommand{\Gr}{\mathrm{Grass}\, (n,2n)}
\newcommand{\eqr}[1]{~\mbox{$(${\rm \ref{#1}}$)$}}

\def\2stack#1#2{\mathrel{\mathop{#1}\limits_{#2}}}
\def\3stack#1#2#3{\mathrel{\mathop{\mathop{#1}\limits_{#2}}\limits_{#3}}}

\title{The Multiplicative Inverse Eigenvalue Problem
over an Algebraically Closed Field}
\date{\today}
\author{Joachim Rosenthal\thanks{Supported 
   in part by NSF grant DMS-96-10389.}\\
  {\small Department of Mathematics}\\
  {\small University of Notre Dame}\\
  {\small Notre Dame, Indiana 46556-5683}\\
  {\small {\em e-mail:} Rosenthal.1@nd.edu} \and
  Xiaochang Wang\\
  {\small Department of Mathematics}\\
  {\small Texas Tech University}\\
  {\small Lubbock, Texas 79409-1042}\\
  {\small {\em e-mail:} mdxia@ttacs.ttu.edu}}

\begin{document} \maketitle
\begin{abstract}
  Let $M$ be an $n\times n$ square matrix and let $p(\lambda)$ be
  a monic polynomial of degree $n$. Let $\cZ$ be a set of
  $n\times n$ matrices. The multiplicative inverse eigenvalue
  problem asks for the construction of a matrix $Z\in\cZ$ such
  that the product matrix $MZ$ has characteristic polynomial
  $p(\lambda)$. 

  In this paper we provide new necessary and sufficient
  conditions when $\cZ$ is an affine variety over an
  algebraically closed field.\bigskip

\noindent
{\bf Keywords:} Eigenvalue completion, inverse eigenvalue
problems, dominant morphism theorem.
\end{abstract}

\Section{Introduction}

Inverse eigenvalue problems involving partially specified
matrices have drawn the attention of many researchers. The
problems are of significance both from a theoretical point of
view and from an applications point of view. For background
material we refer to the monograph by Gohberg, Kaashoek and van
Schagen~\cite{go95}, the recent one by Xu~\cite{xu98} and the
survey article by Chu~\cite{ch98}.

The multiplicative eigenvalue problem asks for conditions which
guarantee that the spectrum of a certain matrix $M$ can be made
arbitrarily through pre-multiplication by a matrix from a certain
set. To be precise let $\F$ be an arbitrary field.  Let $\mat$ be
the space of all $n\times n$ matrices defined over the field
$\F$. We will identify $\mat$ with the vector space $\F^{n^2}$.
Let $\cZ \subset \mat$ be an arbitrary subset and let $M\in\mat$
be a fixed matrix. Then the (right) multiplicative inverse
eigenvalue problem in its general form asks:

\begin{prob} \label{problem}
  Given a monic polynomial $p(\lambda)$ of degree $n$. Is there a
  $n\times n$ matrix $Z\in \cZ$ such that $MZ$ has characteristic
  polynomial
  $$
  \det (\lambda I-MZ)=p(\lambda).
  $$
\end{prob}

The formulation of the left multiplicative inverse eigenvalue
problem is analogous, seeking a matrix $Z\in \cZ$ such that $ZM$
has characteristic polynomial $p(\lambda)$. The left and the
right multiplicative inverse eigenvalue problem are equivalent to
each other because of the identity:
$$
\det (\lambda I-ZA)=\det (\lambda I-A^tZ^t).
$$

In its general form Problem~\ref{problem} is an `open end
problem' and until this point only very particular situations are
well understood. E.g. we would like to mention the well known
result by Friedland~\cite{fr77a} who considered the set
$\cZ=\mathcal{D}$ of diagonal matrices. Friedland did show in
this case by topological methods that Problem~\ref{problem} has
an affirmative answer if the base field $\F$ consists of the
complex numbers $\C$. This diagonal perturbation result was later
generalized by Dias da~Silva~\cite{di86} to situations where the
base field can be any algebraically closed field.

The result which we are going to derive in this paper can be
viewed as a large generalization of Friedland's result.
Specifically we will deal with the situation where
$\cZ\subset\mat$ represents an arbitrary affine variety over an
arbitrary algebraically closed field $\F$. Under these
assumptions we will derive necessary and sufficient conditions
(Theorem~\ref{main}) which will guarantee that
Problem~\ref{problem} has a positive answer for a `generic set'
of matrices $M$ and a `generic set' of monic polynomials
$p(\lambda)$ of degree~$n$.

The techniques which we use in this paper have been developed by
the authors in the context of the additive inverse eigenvalue
problem~\cite{by93,he97,ro00a} and in the context of the pole
placement problem~\cite{ro97}. 

The major tool from algebraic geometry which we will use is the
`Dominant Morphism Theorem' (see Theorem~\ref{dominant}). This
powerful theorem necessitates that the base field is
algebraically closed. The situation over a non-algebraically
closed field seems to be much more complicated. Some new
techniques applicable over the real numbers have been recently
reported by Drew, Johnson, Olesky and van den
Driessche~\cite{dr00}.

\Section{Preliminaries}
For the convenience of the reader we provide a summary of results
which will be needed to establish the new results of this paper.

Denote by $\sigma_i(M)$ the $i$-th elementary symmetric function
in the eigenvalues of $M$, i.e.  $\sigma_i(M)$ denotes up to sign
the $i$-th coefficient of the characteristic polynomial of $M$.
Crucial for our purposes will be the {\em eigenvalue assignment
  map}
\begin{equation}
\psi :\ \cZ \longrightarrow \F^n,\hspace{3mm} Z \longmapsto
(-\sigma_1(MZ),\dots,(-1)^n\sigma_n(MZ)).
\end{equation}
$\psi$ is a morphism in the sense of algebraic geometry.  By
identifying a monic polynomial $\lambda^n+b_1\lambda^{n-1}+\cdots
+b_n$ with the point $(b_1,\dots b_n)\in \F^n$ we can also write
\begin{equation}
\psi (Z)=\det (\lambda I-MZ).
\end{equation}

Crucial for the proof of the main result (Theorem~\ref{main})
will be the Dominant Morphism Theorem. The following version can
be immediately deduced from~\cite[Chapter AG, \S 17, Theorem
17.3]{bo91b1}.

\begin{pr}                           \label{dominant}
  Let $f:\cZ \rightarrow {\mathcal Y}$ be a morphism of affine
  varieties over an algebraically closed field.  Then the image
  of $f$ contains a nonempty Zariski open set of ${\mathcal Y}$
  if and only if the Jacobian $df_{_Z}:{T}_{_{Z,\cZ}}\rightarrow
  {T}_{_{f(Z),{\mathcal Y}}}$ is onto at some smooth point $Z$ of
  $\cZ$, where ${T}_{_{X,{\mathcal X}}}$ is the tangent space of
  ${\mathcal X}$ at the point $X$.
\end{pr}

There are classical formulas, sometimes referred to as Newton
formulas, which express the elementary symmetric functions
$\sigma_i(M)$ uniquely as a polynomial in the power sum symmetric
functions
$$
p_i:=\lambda_1^i+\cdots+\lambda_n^i=\tr(M)^i.
$$
To be precise one has the formula (see e.g.~\cite{ma79})
$$
\sigma_i(M)=\frac{1}{n!}  \left(\begin{array}{ccccc}
    p_1  &   1  &   0  & \ldots &     0  \\
    p_2  &  p_1 &   2  &        & \vdots \\
    \vdots&\ddots&\ddots&\ddots  & \vdots \\
    \vdots&      &\ddots& p_1  & n-1  \\
    p_n &\ldots&\ldots& p_2 & p_1
\end{array}\right),
$$
which induces an isomorphism
$\F^n\rightarrow\F^n,(p_1,\ldots,p_n)\mapsto(\sigma_1,\ldots,\sigma_n)$.
Based on this we equally well can study the map
\begin{equation}
\phi :\ \cZ \longrightarrow\F^n,\hspace{3mm}M
\longmapsto \left({ \tr(MZ),\ldots,\tr((MZ)^n) }\right).
\end{equation}

We will use the following result from~\cite{he97}:
\begin{pr}                           \label{mainpr}
  Let ${\mathcal L}\subset \mat$ be a linear sub-space of
  dimension $\geq n$, ${\mathcal L}\not\subset sl_n$ (i.e.
  ${\mathcal L}$ contains an element with nonzero trace).  Define
  $$
  \pi (M)=(m_{11},m_{22},\dots,m_{nn})
  $$
  the projection onto the diagonal entries.  Then there exists
  a $S\in Gl_n$ such that
  $$ 
  \pi(S{\mathcal L}S^{-1})=\F^n.
  $$
\end{pr}

It is possible to `compactify' the problem.  For this consider
the identity
\begin{equation}
  \label{ident}
  \det (\lambda I-MZ)=\det \left[\begin{array}{cc}I& Z\\
  M& \lambda I
\end{array}\right].
\end{equation}
Denote by $\mathrm{Grass}\, (k,n)$ the Grassmann manifold
consisting of all $k$-dimensional linear subspaces of $\F^n$.
Algebraically $\mathrm{Grass}\, (k,n)$ has the structure of a
smooth projective variety. In the sequel we will identify
$\mathrm{rowsp}[I\ Z]$ with a point in the Grassmannian $\Gr$. By
identifying $\mathrm{rowsp}[I\ Z]$ with $Z\in\mat$, we can say
that $\cZ\subset \Gr$. Let $\bar{\cZ}$ be the projective closure
of $\cZ$ in $\Gr$. Every element in $\bar{\cZ}$ can be simply
represented by a subspace of the form $\mathrm{rowsp}[Z_1\ Z_2]$,
where the $n\times n$ matrix $Z_1$ is not necessarily invertible.
$\mathrm{rowsp}[Z_1\ Z_2]$ describes an element of $\cZ$ if and
only if $Z_1$ is invertible. For any element $\mathrm{rowsp}[Z_1\ 
Z_2]\in \bar{\cZ}$, define $\bar{\psi}: \bar{\cZ} \longrightarrow
\mathbb{P}^n$
\begin{equation} \label{barmap}
\bar{\psi}([Z_1\ Z_2])=\det \left[\begin{array}{cc}Z_1& Z_2\\
  M& \lambda I
\end{array}\right].
\end{equation}
where a polynomial $b_0\lambda^n+ b_1\lambda^{n-1}+\cdots +b_n$
is identified with the point $(b_0,b_1,\dots,b_n)\in
\mathbb{P}^n$. Recall that the Pl\"{u}cker coordinates of
$\mathrm{rowsp}[Z_1\ Z_2]\in \Gr$ are given by the full size
minors $[Z_1\ Z_2]$, and by considering the Pl\"{u}cker
coordinates as the homogeneous coordinates of points in
$\mathbb{P}^N$, $N={2n\choose n}-1$, one has an embedding
$\Gr\subset \mathbb{P}^N$ which is called Pl\"{u}cker embedding.
Under the Pl\"{u}cker coordinates, (\ref{barmap}) becomes
\begin{equation} \label{barmap2}
\bar{\psi}([Z_1\ Z_2])=\det \left[\begin{array}{cc}Z_1& Z_2\\
  M& \lambda I
\end{array}\right]=\sum_{i=0}^N z_i m_i(\lambda)
\end{equation}
where $\{z_i\}$ are $n\times n$ minors of $[Z_1\ Z_2]$ and
$m_i(\lambda)$ is the cofactor of the $z_i$ in the determinate
of\eqr{barmap}. $\bar{\psi}$ is undefined on the elements where
$$
\det \left[\begin{array}{cc}Z_1& Z_2\\
    M& \lambda I
\end{array}\right]=0.
$$
So $\bar{\psi}$ is a rational map.

\Section{New Results}
The next theorem constitutes the main result of this paper. As
stated in the introduction we will identify the set $\mat$ with
the vector space $\F^{n^2}$ and we will identify the set of monic
polynomials of degree $n$
$$
\lambda^n+b_1\lambda^{n-1}+\cdots +b_n
$$
with the vector space $\F^n$. If $V$ is an arbitrary
$\F$-vector space one says that $U\subset V$ forms a generic set
if $U$ contains a non-empty Zariski open subset. Over the complex
or real numbers a generic set is necessarily dense with respect
to the natural topology. The Dominant Morphism
Theorem~\ref{dominant} states that the image of an algebraic
morphism forms a generic set as soon as the linearization around
a smooth point is surjective and if the field is algebraically
closed.

If Problem~\ref{problem} has a positive answer for a generic set
of matrices inside $\mat$ and a generic set of monic polynomials
then we will say that Problem~\ref{problem} is generically
solvable. With this preliminary we have the main result of this
paper:

\begin{theorem}                    \label{main}
  Let $\cZ\subset\mat$ be an affine variety over an algebraically
  closed field~$\F$. Then Problem~\ref{problem} is generically
  solvable if and only if $\dim \cZ\geq n$ and $\det (Z)$ is not
  a constant function on $\cZ$.
\end{theorem}

\begin{proof}
  The conditions are obviously necessary.  So we only need to
  prove the sufficiency. Assume that $\dim \cZ\geq n$ and $\det
  (Z)$ is not a constant on $\cZ$.  Then there exists a curve
  $Z(t)\subset \cZ$ such that
  $$
  \frac{d}{dt} \det Z(t)|_{t=0}\neq 0,
  $$
  $Z(0)=Z_0$ is a smooth point of $\cZ$, and $\det Z_0\neq 0$.
  
  Let $Z(t)=Z_0+tL+O(t^2)$ where $L\in {T}_{_{Z_0,\cZ}}$. Then
  $$
  \det Z(t)=\det Z_0\det (I+tZ_0^{-1}L+O(t^2))= \det
  Z_0(1+t\tr Z_0^{-1}L+O(t^2))
  $$
  and
  $$
  \frac{d}{dt} \det Z(t)|_{t=0}=\det Z_0\tr Z_0^{-1}L \neq 0,
  $$
  i.e.
  $$
  Z_0^{-1}{T}_{_{Z_0,\cZ}}\not\subset sl_n.
  $$
  By Proposition~\ref{mainpr}, there exists a $S\in Gl_n$ such
  that
  $$
  \pi(SZ_0^{-1}{T}_{_{Z_0,\cZ}}S^{-1})=\F^n.
  $$
  Let
  \begin{equation}   \label{d}
  D:=\left[
  \begin{array}{cccc}
  1&&&\\
  &2&&\\
  &&\ddots&\\
  &&&n
  \end{array}\right]
  \end{equation}
  and
  $$
  M:=S^{-1}DSZ_0^{-1}.
  $$
  Then for any curve through $Z_0$
  $$
  Z(t)=Z_0+tL+O(t^2)\subset \cZ,\ \ \ L\in {T}_{_{Z_0,\cZ}},
  $$
  we have
  \begin{eqnarray*}
  \lim_{t\rightarrow 0}\frac{ \tr (MZ(t))^i-\tr (MZ_0)^i}{t}
  &=&\lim_{t\rightarrow 0}
  \frac{ \tr (MZ_0+tML+O(t^2))^i-\tr (MZ_0)^i}{t}\\
  &=&i\cdot \tr ((MZ_0)^{i-1}ML)\\
  &=&i\cdot \tr (D^iSZ_0^{-1}LS^{-1}).
  \end{eqnarray*}
  
  Let
\begin{equation} \label{v}
V=D\left[\begin{array}{cccc}
1&1&\cdots&1\\
1&2&\cdots&2^{n-1}\\
\vdots&\vdots&  &\vdots\\
1&n&\cdots&n^{n-1}
\end{array}\right]D.
\end{equation}
Then $V$ is invertible and the Jacobian
$d\phi_{_{Z_0}}:{T}_{_{Z_0,\cZ}}\mapsto \F^n$
\begin{eqnarray*}
d\phi_{_{Z_0}}(L)
&=&(\tr(DSZ_0^{-1}LS^{-1}),
2\tr(D^2SZ_0^{-1}LS^{-1}),\dots,n\tr(D^nSZ_0^{-1}LS^{-1}))\\
&=&\pi(SZ_0^{-1}LS^{-1})V
\end{eqnarray*}
is onto.  By the Dominant Morphism Theorem~\ref{dominant},
$\phi(\cZ)$ contains a nonempty Zariski open set of $\F^n$, so
does $\psi(\cZ)$.

Since the set of $M$'s such that $\psi$ is almost onto is a
Zariski open set, and we just showed that it is nonempty, $\psi$
is almost onto for a generic set of matrices $M$.
\end{proof}

Next we consider the number of solutions of Problem~\ref{problem}
when $\dim \cZ=n$. For this we introduce an important technical
concept.

\begin{de}
  A matrix $M$ is called $\cZ$-nondegenerate for the right
  multiplicative inverse eigenvalue problem if
  \begin{equation}                       \label{non-deg}
  \det \left[\begin{array}{cc}Z_1& Z_2\\
      M& \lambda I
  \end{array}\right]\neq 0
  \end{equation}
  for any $\mathrm{rowsp}[Z_1, Z_2]\in \bar{\cZ}\subset \Gr$.
\end{de}
So if $M$ is $\cZ$-nondegenerate, then the map $\bar{\psi}$
defined by\eqr{barmap} becomes a morphism. In this situation we
can say even quite a bit more:

\begin{theorem}                    \label{main2}
  If $M$ is $\cZ$-nondegenerate, and $\dim \cZ=n$, then
  Problem~\ref{problem} is solvable for any monic polynomial
  $p(\lambda)$ of degree $n$. Moreover, when counted with
  multiplicities, the number of matrices inside $\cZ$ which
  results in a characteristic polynomial $p(\lambda)$ is exactly
  equal to the degree of the projective variety $\bar{\cZ}\subset
  \Gr$ when viewed under the Pl\"{u}cker embedding $\Gr \subset
  \mathbb{P}^N$.
\end{theorem}

\begin{proof}
  We will repeatedly use the projective 
  dimension theorem~\cite[Charter I, Theorem 7.2]{ha77} which says
  that if $X$ and $Y$ are $r$-dimensional and $s$-codimensional
  projective varieties, respectively, then $\dim X\cap Y\geq r-s$.
  In particular, $X\cap Y$ is not empty if $r\geq s$.

  Let
  $$
  K=\{(z_0,\dots,z_{_N})\in \mathbb{P}^N|\sum_{i=0}^Nz_i
  m_i(\lambda)=0\}.
  $$
  Then $K$ must have co-dimension $n+1$ because of the
  condition $K\cap \cZ=\emptyset$.  Therefore the linear equation
  \begin{equation}\label{linear}
   \sum_{i=0}^N z_i m_i(\lambda)=p(\lambda).
  \end{equation}
  has solutions in $\mathbb{P}^N$ for any $p(\lambda)\in
  \mathbb{P}^n$, and the set of all solutions for each 
  $p(\lambda)$ is in the form of $z_p+K$ 
  where $z_p$ is a particular solution;  i.e. the solution set is 
  given by $K_p-K$ where $K_p$ is the unique $n$-codimensional 
  projective subspace through $z_p$ and $K$.
  Since
  $K\cap \bar{\cZ}=\emptyset$, we must have
  $$
  \dim K_p\cap \bar{\cZ}=0
  $$
  and by B\'ezout's Theorem~\cite{vo84}, there are $\deg
  \bar{\cZ}$ many points in $K_p\cap \bar{\cZ}$ counted with
  multiplicities.  If $p(\lambda)$ is a monic polynomial of degree
  $n$, then
  from\eqr{barmap} one can see that all the solutions are in
  $\cZ$.
\end{proof}

An immediate application of Theorem~\ref{main2} is a result of
Friedland~\cite{fr77a}: Let $\cZ$ be the set of all diagonal
matrices. Then closure $\bar{\cZ}$ of $\cZ$ inside the Grassmann
variety $\Gr$ is isomorphic to the product of $n$ projective
lines:
$$
\mathbb{P}^1\times\cdots \times \mathbb{P}^1
$$
As shown in~\cite{by93} the degree of $\bar{\cZ}$ is then
equal to $n!$. Moreover all points of $\bar{\cZ}$ are of the form
$\mathrm{rowsp}[Z_1\ Z_2]$ where $Z_1$ and $Z_2$ are given by
$$
Z_1=\left[\begin{array}{cccc}
    z_{11}& 0    &\cdots & 0     \\
    0     &z_{12}&\cdots & 0     \\
    \vdots&\vdots&\ddots & \vdots\\
    0 & 0 &\cdots & z_{1n}
\end{array}\right],\
Z_2=\left[\begin{array}{cccc}
    z_{21}& 0    &\cdots & 0     \\
    0     &z_{22}&\cdots & 0     \\
    \vdots&\vdots&\ddots & \vdots\\
    0 & 0 &\cdots & z_{2n}
\end{array}\right].
$$
In these matrices, $(z_{1i},z_{2i})$ represent the homogeneous
coordinates of the $i$th projective line $\mathbb{P}^1$.

In order to apply Theorem~\ref{main2} we have to find the
algebraic conditions which guarantee that a particular matrix $M$
is $\cZ$-nondegenerate, i.e. condition\eqr{non-deg} has to be
satisfied for every element $[Z_1\ Z_2]\in \bar{\cZ}$. For this
let $I$ be a subset of $\{1,2,\dots,n\}$, $J$ be the complement
of $I$, and $|J|$ be the number of elements in $J$. For any point
$[Z_1\ Z_2]\in\bar{\cZ}$, assume
$$
\begin{array}{ll}
  z_{1i}=0&\mbox{for $i\in I$},\\
  z_{1j}\neq 0&\mbox{for $j\in J$}.
\end{array}
$$
Without loss of generality we can take
$$
\begin{array}{ll}
  z_{2i}=1&\mbox{for $i\in I$},\\
  z_{1j}=1&\mbox{for $j\in J$},
\end{array}
$$
and\eqr{barmap} becomes
$$
\bar{\psi}([Z_1\ Z_2])=\pm M_{_I}\lambda^{|J|}+\mbox{lower
  power terms},
$$
where $M_{_I}$ is the principal minor of $M$ consisting of the
$i$th rows and columns, $i\in I$. Furthermore if we take
$$
z_{2j}=0\mbox{ for $j\in J$}
$$
then
$$
\bar{\psi}([Z_1\ Z_2])=\pm M_{_I}\lambda^{|J|}.
$$
Therefore $M$ is $\cZ$-nondegenerate if and only if all the
principal minors of $M$ are nonzero. So we have Friedland's
result~\cite[Theorem 2.3]{fr77a} formulated for an algebraically
closed field: If all the principal minors of $M$ are nonzero,
then the multiplicative inverse eigenvalue problem with
perturbation from the set of diagonal matrices is solvable for
any monic polynomial $p(\lambda)$ of degree $n$, and there are
$n!$ solutions, when counted with multiplicities.


\end{document}